\documentclass[12pt]{article}

\usepackage{amsmath, amssymb}
\usepackage{amsfonts, mathrsfs}
\usepackage[cp1251]{inputenc}
\usepackage[active]{srcltx}
\usepackage[ukrainian, russian, english]{babel}

\usepackage{amsfonts, amssymb, mathrsfs}

\newtheorem{theorem}{Теорема}
\newtheorem{lemma}{Лемма}
\newtheorem{corollary}{Следствие}
\newtheorem{proposition}{Предложение}

\begin{document}

\selectlanguage{russian}

\setcounter{page}{1}

\noindent {\small УДК 517.518.83} \vskip 3mm

\noindent \textbf{А.\,С. Сердюк, И.\,В. Соколенко}  {\small (Институт математики НАН Украины, г.~Киев)} \vskip 5mm

\noindent \textbf{A.S. Serdyuk, I.V. Sokolenko} {\small (Institute of Mathematics
of The National Academy of Sciences of Ukraine, Kiev)}\vskip 5mm

\noindent  \textbf{\Large Асимптотические равенства для наилучших приближений  классов бесконечно дифференцируемых функций,  задающихся с помощью модуля непрерывности}
\vskip 10mm

\noindent  \textbf{\Large Asymptotic behavior of best approximations of
classes of infinitely differentiable functions defined by
moduli of continuity}\vskip 10mm

{\small \it \noindent
  Получены  асимптотические оценки для наилучших   приближений тригонометрическими полиномами в метрике пространства $C$ $\ (L_p)\ $  классов  периодических функций, представимых в виде  сверток ядер $\Psi_\beta$,  с  коэффициентами Фурье, убывающими к нулю
   быстрее любой степенной последовательности, с функциями $\varphi\in C\ (\varphi\in L_p)$, модули непрерывности которых не превышают заданной мажоранты $\omega(t)$. Доказано, что в пространствах $C$ и $L_1$ для выпуклых модулей непрерывности $\omega(t)$ полученные оценки являются асимптотически точными.
\hfill}
\vskip 5mm

{\small \it \noindent
We obtain asymptotic estimates for the best approximations by trigonometric polynomials in the metric space $C$ $\ (L_p)\ $ of classes of periodic functions that can be represented as a convolution of kernels $\Psi_\beta$, which Fourier coefficients tend to zero faster than any power sequence, with functions $\varphi\in C\ (\varphi\in L_p),$ which moduli of continuity do not exceed a fixed majorant $\omega(t)$. It is proved that in the spaces $C$ and $L_1$  the obtained estimates are asymptotically exact for convex moduli of continuity $\omega(t)$.
\hfill}
\vskip 1.5mm

\begin{center}
\textbf{Постановка задачи и формулировка основных результатов}
\end{center}

Пусть $C$ и $\ L_p,\ 1\le p\le\infty,$~--- пространства $2\pi$-периодических  функций
 со стандартными  нормами $\|\cdot\|_{C}$  и $\|\cdot\|_{L_p}$.
Пусть, далее,    $L_{\beta }^{\psi }$ ---  множество функций  $f\in L_1$, представимых  сверткой
\begin{equation}\label{1}
f(x)=\frac {a_0(f)}2 + \frac 1{\pi }\int\limits_{-\pi }^{\pi }
\varphi(x-t)\Psi _{\beta }(t)dt,\ \ \ \varphi\in L_1,\ \ \ \varphi\bot1,
\end{equation}
с суммируемым ядром $\Psi_\beta(t),$ ряд Фурье которого имеет вид
$\
\Psi _{\beta }(t)\sim \sum\limits_{k=1}^{\infty}\psi (k)\cos \left(kt -\frac {\beta\pi }2\right), \   \psi(k)>0, \   \beta\in\mathbb R. \
$
Равенство (\ref{1}) понимается как равенство двух функций из $L_1$, т.е. почти для всех $x\in \mathbb R.$
Согласно   Степанцу \cite[\S 3.7]{Stepanets2002_1}, функцию $\varphi $  в равенстве $(\ref{1})$  называют
$(\psi,\beta)$-про\-из\-водной  функции $f$ и обозначают
$f_{\beta}^{\psi}$ \ \ ($\varphi(x)=f_{\beta}^{\psi}(x)$).
 Если $f\in L_{\beta }^{\psi }$  и, кроме того, $f_{\beta
}^{\psi }\in {\mathfrak N},$ где ${\mathfrak N} \subseteq L_1,$ то
записывают $f\in L_{\beta }^{\psi }{\mathfrak N}$.
 Также полагают \ $C_{\beta}^{\psi}=L_{\beta}^{\psi}\cap
C$\ \  и \ \ $C_{\beta}^{\psi}\mathfrak N=L_{\beta}^{\psi}\mathfrak N\cap
C$.

Пусть $\omega(t)$ --- фиксированный модуль непрерывности, т.е. непрерывная неубывающая и полуаддитивная  при всех $t\ge0$ функция, в нуле равная нулю.
В   качестве  $\mathfrak N$ рассматриваются множества
$
H_{\omega_X}=\{\varphi\in X: \omega_X(\varphi;t)\le \omega(t), \  t\ge 0\},
$
где
$\omega_X(\varphi;t)=\sup\limits_{|h|\le t}\|\varphi(\cdot+h)-\varphi(\cdot)\|_X$
--- модуль непрерывности функции $\varphi$ в пространстве $X,$ а $X$ есть $C$ либо $L_p, 1{\le} p{<}\infty.$
При $\psi(k){=}k^{-r} $ и $ \beta{=}r,  r{\in} \mathbb N,$  классы $C_{\beta}^{\psi} H_{\omega_C}$ совпадают с известными классами  $W^r H_{\omega_C}$ $2\pi$-периодических $r$ раз непрерывно дифференцируемых функций, $r$-е производные которых принадлежат классу $H_{\omega_C}.$

Не уменьшая общности, последовательность $\psi(k)$, определяющую
классы $L_{\beta}^{\psi}\mathfrak N$ и $C_{\beta}^{\psi}\mathfrak N$, можно считать следом на  $\mathbb{N}$ некоторой непрерывной функции  $\psi(t)$ непрерывного аргумента $t\in[1,\infty)$. Множество всех положительных непрерывных выпуклых вниз функций $\psi(t),\ t\in[1,\infty)$, для которых $\mathop{\lim}\limits_{t\rightarrow\infty}\psi(t)=0,\ $ обозначим через $\mathfrak M$.
 Каждой функции
$\psi \in \mathfrak M$ поставим  в соответствие  характеристики (см. \cite[\S 3.12]{Stepanets2002_1})
$\eta(t)=\eta(\psi;t)=\psi^{-1}\left( {\psi(t)}/2\right)$ и $ \mu(t)=\mu(\psi;t)= {t}/(\eta(t)-t),$
где $\psi^{-1}$~--- функция, обратная к $\psi$, и выделим подмножество
$\mathfrak M_{\infty}^+$  функций
 $\psi $ из $ \mathfrak M$,   для каждой из которых
$\mu(t)=\mu(\psi;t)$ монотонно и неограниченно возрастает:
$ \mu(\psi;t)\uparrow\infty.$
Естественным  представителем множества $\mathfrak M_\infty^+$ является функция
\begin{equation}\label{13}
\psi(t)=e^{-\alpha t^{r}}, \ \ \ \alpha>0, \ \ \ r>0,
\end{equation}
для которой  характеристика $\mu(t)$ имеет вид
\begin{equation}\label{15}
\mu(t)=\mu(\psi;t)=\left(\left(1+\frac{\ln 2}{\alpha
t^{r}}\right)^{1/r}-1 \right)^{-1}.
\end{equation}
Ядро $\Psi_\beta(t)$, порождаемое функцией $\psi(t)$ вида $(\ref{13})$,
 называют обобщенным ядром Пуассона и обозначают $P_{\beta}^{\alpha,r}(t),$ т.е.
\begin{equation}\label{14}
    P_{\beta}^{\alpha,r}(t)=\sum _{k=1}^{\infty}e^{-\alpha k^r}\cos \left(kt -\frac {\beta\pi }2\right),\ \ \ \alpha>0, \ \ \ r>0, \ \ \ \beta\in\mathbb R.
\end{equation}
В этом случае  классы  $C_{\beta}^{\psi}\mathfrak N$ и
$L_{\beta}^{\psi}\mathfrak N$ обозначают соответственно через $C_{\beta}^{\alpha,r}\mathfrak N$ и $L_{\beta}^{\alpha,r}\mathfrak N$.
Изучению асимптотического поведения уклонений сумм Фурье от функций из классов $C_{\beta}^{\alpha,r}\mathfrak N$ и $L_{\beta}^{\alpha,r}\mathfrak N$ в метриках пространств $C$ и $L_p$ посвящены работы   Степанца \cite{Stepanets1984,Stepanets2001} и   Теляковского \cite{Tel1989}. Исследованию асимптотического поведения уклонений обобщенных операторов Абеля--Пуассона, порождаемых ядрами вида (\ref{14}) при $\beta{=}0$, от функций из классов Липшица ${\rm Lip}_1\gamma, 0{<}\gamma{<}r{\le}1,$ посвящена  работа  Фалалеева \cite{Falaleev2000}.

Как следует из предложение 5.14.1 работы \cite{Stepanets2002_1}, если $\psi\in\mathfrak M_\infty^+$ и  $f\in C^\psi_\beta,\ \beta\in\mathbb R,$  то   функция  $f$ является бесконечно дифференцируемой. С другой стороны,  из     \cite[теорема 2]{Stepanets_Serdyuk_Sh_2008} вытекает, что если $2\pi$-периодическая функция $f$ является бесконечно дифференцируемой, то существует  $\psi\in\mathfrak M_\infty^+$ для которой $f$ принадлежит $ C^\psi_\beta$ при всех действительных $\beta$. Иными словами, если $D^\infty$ --- множество $2\pi$-периодических бесконечно дифференцируемых функций, то $\ \bigcup\limits_{\psi\in\mathfrak M^+_\infty} C^\psi_\beta=D^\infty \ $ при произвольном действительном $\beta.$

Обозначим через  $\ {\cal T}_{2n-1}\ $   подпространство всех тригонометрических полиномов
$T_{n{-}1}(x){=}\sum\limits_{k=0}^{n{-}1} (\alpha_k\cos kx{+}\beta_k\sin kx),   \alpha_k,\beta_k{\in} \mathbb{R}, $
порядка не выше $n{-}1$. Данная работа посвящена нахождению асимптотических при $n{\rightarrow}\infty$ оценок наилучших приближений
 \begin{equation*}\label{7}
 E_n( C_{\beta}^\psi H_{\omega_C})_C=\sup\limits_{f\in C_{\beta}^\psi H_{\omega_C}} E_n(f)_C= \sup\limits_{f\in C_{\beta}^\psi H_{\omega_C}}\inf\limits_{T_{n-1}\in{\cal T}_{2n-1}}
 \|f-T_{n-1}\|_C,
\end{equation*}
и
\begin{equation*}\label{7'}
 E_n( L_{\beta}^\psi H_{\omega_{L_p}})_{L_p}=\sup\limits_{f\in L_{\beta}^\psi H_{\omega_{L_p}}} E_n(f)_{L_p}= \sup\limits_{f\in L_{\beta}^\psi H_{\omega_{L_p}}}\inf\limits_{T_{n-1}\in{\cal T}_{2n-1}}
 \|f-T_{n-1}\|_{L_p},
\end{equation*}
 в случае, когда $\psi\in\mathfrak M_\infty^+$ и $\beta\in\mathbb R.$

Для функциональных классов, задающихся с помощью
модуля непрерывности, на сегодняшний день
существует относительно небольшое количество работ, в которых
найдены точные  или асимптотические равенства для наилучших приближений таких классов в равномерной  или интегральной
метриках.
   Корнейчук \cite{Korneychuk1970}  вычислил  точные значения величин $E_n(W^rH_{\omega_C})_C$ и
$E_n(W^rH_{\omega_C})_{L_1}, r, n\in \mathbb N,$ для выпуклых модулей
непрерывности $\omega(t)$.
Впоследствии эти результаты были распространены на
классы периодических функций, порождающиеся
не увеличивающим  осцилляцию ядром $K$ ($K\in {\rm CVD}$)
(см. библиографию в \cite{Shevaldin1994}).
  Шевалдин  \cite{Shevaldin1994} для ядер $K$, удовлетворяющих введенному им условию $B_{2n}$, нашел оценки снизу  колмогоровских поперечников классов сверток ядер  $K$ с функциями из $H_{\omega_C}$ в равномерной метрике при условии выпуклости функции $\omega(t)$ и показал, что они в ряде случаев являются точными. Эти результаты дают оценки снизу для наилучших равномерных приближений соответствующих классов. Также в \cite{Shevaldin1994} показано, что ядро Пуассона $P_0^{\alpha,1}$ при $\alpha\ge \ln(1+2\pi)$ удовлетворяет условию $B_{2n} $ для каждого $ n\in\mathbb N,$   но в то же время это ядро    не является CVD-ядром при $\alpha=\ln(1+2\pi)$ ($P_0^{\alpha,1}\notin {\rm CVD}$).

В  работах \cite{Serdyuk_Sokolenko2010_2,Serdyuk_Sokolenko2011_1}
рассмотрен  линейный метод  $U^\ast_{n-1}$ для приближения функций из классов интегралов Пуассона $C_{\beta}^{\alpha,1} H_{\omega_C}$ и $  L_{\beta}^{\alpha,1} H_{\omega_{L_1}}, \alpha{>}0,\ \beta{\in}\mathbb R, $  и решена задача Колмогорова-Никольского, состоящая в нахождении асимптотических при $n{\rightarrow}\infty$ равенств для величин
$
{\cal E}(C_{\beta}^{\alpha,1} H_{\omega_C}; U_{n-1}^\ast)_C=\sup\limits_{f\in C_{\beta}^{\alpha,1} H_{\omega_C}}\|f -U^\ast_{n-1}(f)\|_C
$
 и
 $
 {\cal E}(L_{\beta}^{\alpha,1} H_{\omega_{L_1}}; U_{n-1}^\ast)_{L_1}=\sup\limits_{f\in L_{\beta}^{\alpha,1} H_{\omega_{L_1}}}\|f -U^\ast_{n-1}(f)\|_{L_1}.
 $
На основании найденных оценок   были получены асимптотические равенства для наилучших приближений  $ E_n( C_{\beta}^{\alpha,1} H_{\omega_C})_C$ и $ E_n(L_{\beta}^{\alpha,1} H_{\omega_{L_1}})_{L_1}$ в случае выпуклости модуля непрерывности $\omega(t)$. При этом оказалось, что  величины  ${\cal E}(C_{\beta}^{\alpha,1} H_{\omega_C}; U_{n-1}^\ast)_C,$ $\ {\cal E}(L_{\beta}^{\alpha,1} H_{\omega_{L_1}}; U_{n-1}^\ast)_{L_1}$, $ E_n( C_{\beta}^{\alpha,1} H_{\omega_C})_C$ и $ E_n(L_{\beta}^{\alpha,1} H_{\omega_{L_1}})_{L_1}$ асимптотически равны между собой.
В \cite{Serdyuk_Sokolenko2011_2}   результаты по наилучшим приближениях были распространены   на  классы $C^\psi_\beta H_{\omega_C}$ и  $L^\psi_\beta H_{\omega_{L_1}}$, когда   $\psi(k)$ удовлетворяют условию Д'Aламбера  $\lim\limits_{k\rightarrow\infty}  {\psi(k+1)}/{\psi(k)}=q,\ $ $ 0{<}q{<}1.$ При выполнении этого условия классы $C^\psi_\beta H_{\omega_C}$ состоят из функций, допускающих регулярное продолжение в полосу $|{\rm Im} z|\le \ln (1/q)$.

В данной работе   продолжены исследования работ \cite{Serdyuk_Sokolenko2010_2}--
\cite{Serdyuk_Sokolenko2011_2}. В ней получены   асимптотические  равенства для величин  $\ E_n( C_{\beta}^{\psi}
H_{\omega_C})_C$ и $\ E_n( L_{\beta}^{\psi}
H_{\omega_{L_1}})_{L_1}$ при условии, что  $\psi\in\mathfrak M_\infty^+$, $\beta\in\mathbb R$ и $\omega(t)$ ---  выпуклый модуль непрерывности.  Отметим, что порядковые оценки  этих величин   известны  (см., например, \cite[\S7.4]{Stepanets2002_2}):
$
E_n( C_{\beta}^\psi H_{\omega_C})_{C}\asymp E_n( L_{\beta}^\psi H_{\omega_{L_1}})_{L_1}\asymp \psi(n)\omega\left(  1/n\right).
$

  Основной результат работы составляет  следующее утверждение.

\begin{theorem}\label{theorem1}
Пусть $\psi\in \mathfrak M_{\infty}^+,$  $\ \beta
\in\mathbb{R}\ $   и $\ \omega(t)$~---  произвольный модуль непрерывности.
Тогда при  $\mu(n)>2$
\begin{equation}\label{8}
E_n( C_{\beta}^\psi H_{\omega_C})_{C} = \frac{2}{\pi}\theta_n^{(1)}\psi(n)e_n(\omega) +O(1)\gamma_n(\psi)\omega\left(\frac 1n\right),
\end{equation}
\begin{equation}\label{9}
E_n( L_{\beta}^\psi H_{\omega_{L_1}})_{L_1} = \frac{2}{\pi}\theta_n^{(2)}\psi(n)e_n(\omega) +O(1)\gamma_n(\psi)\omega\left(\frac 1n\right),
\end{equation}
\begin{equation}\label{10}
E_n( L_{\beta}^\psi H_{\omega_{L_p}})_{L_p}\le \frac{2}{\pi}\psi(n)e_n(\omega) +O(1)\gamma_n(\psi)\omega\left(\frac 1n\right),\ \ \ 1<p<\infty,
\end{equation}
где
\begin{equation}\label{11}
    e_n(\omega)=\int\limits_0^{\pi/2} \omega \bigg( \frac{2t}{n} \bigg) \sin t dt,
\end{equation}
\begin{equation}\label{12}
    \gamma_n(\psi)=  \psi(n+1)\left(\frac1n+\frac1{\mu(n)-2}\right)+\psi(3n)\left(1+\ln^+(\eta(n)-n)\right),
\end{equation}
$\ln^+t=\max\{0;\ln t\},$ величины  $ \theta_n^{(i)}=\theta_n^{(i)}(\omega),$ $ i=1,2,$ таковы, что $\theta_n^{(i)}=1$  для выпуклых модулей непрерывности $\omega(t)\ $   и $ 2/3\le \theta_n^{(1)} \le1,$  $ 1/2\le \theta_n^{(2)} \le1$ для произвольных $\ \omega(t)$. В формулах (\ref{8})--(\ref{10})  и всюду далее
     $O(1)$ --- величины равномерно ограниченные по всем рассматриваемым параметрам.
\end{theorem}

Приведем некоторые замечания и следствия из теоремы 1. В случае выпуклых модулей непрерывности $\omega(t)$ формулы (\ref{8}) и (\ref{9}) являются асимптотическими равенствами тогда и только тогда, когда  выполняется условие
\begin{equation}\label{12'}
    \frac{\psi(3n)}{\psi(n)}=o\left(\frac1{1+\ln^+(\eta(n)-n)}\right).
\end{equation}
Из \cite[с.~143]{Stepanets2002_2} известно, что
$\ \psi(3n)/\psi(n)=o(1)$ при $\psi\in\mathfrak M_\infty^+.$ Поэтому, если $\psi(t)$ такова, что
$    \eta(\psi;t)-t\le K<\infty,\ \ \ t\in[1,\infty),\ $  то условие (\ref{12'}) для таких $\psi(t)$  выполняется всегда. Если же разность $\eta(\psi;t)-t $ не ограничена, то соотношение (\ref{12'}) может и не выполнятся.

Для оценки величины $\gamma_n(\psi)$ вида (\ref{12}) может быть полезным следующее утверждение.

\begin{proposition}
Пусть заданные на $[1,\infty)$ функции $\psi_1(t)$ и $\psi_2(t)$    положительны, непрерывны и строго убывают к нулю, а $\varphi(t)= {\psi_2(t)}/{\psi_1(t)}$ монотонно неубывает. Тогда
$ \   \eta(\psi_1;t)\le \eta(\psi_2;t), \ \   t\ge1.$
\end{proposition}

Доказательство вытекает из  соотношений:
\begin{equation*}\label{1l3}
    \psi_1(\eta(\psi_1;t))=\frac12\psi_1(t)=\frac12\frac{\psi_2(t)}{\varphi(t)}\ge\frac{\frac12\psi_2(t)}{\varphi(\eta(\psi_2;t))}=
    \frac{\psi_2(\eta(\psi_2;t))}{\varphi(\eta(\psi_2;t))}=\psi_1(\eta(\psi_2;t)).
\end{equation*}

Наряду с характеристиками $\eta(t)$ и $\mu(t)$ для функций $\psi$ из множества  $\mathfrak M$ рассматривают также (см. \cite[\S3.12]{Stepanets2002_1})   характеристику $\alpha(t)=\alpha(\psi;t)=\frac{\psi(t)}{t|\psi'(t)|}$ $(\psi'(t){\mathop{=}\limits^{\rm df}}\psi'(t{+}0))$. С её помощью  выделим из $\ \mathfrak M\ $ подмножества
$\
S^+=\left\{\psi\in\mathfrak M:\ \frac1{\alpha(t)}\nearrow\right\}\ \ $ и $\ S^+_\infty=\left\{\psi\in\mathfrak M:\ \frac1{\alpha(t)}\uparrow\infty\right\}.
$
Из теоремы 3.12.1 работы \cite{Stepanets2002_1} вытекает, что $S_\infty^+\subset \mathfrak M_\infty^+.$
Ниже формулируется  утверждение, полезное для проверки условия принадлежности множеству $S_\infty^+$, а значит, и $ \mathfrak M_\infty^+$.

\begin{proposition}
   Пусть  $\psi_1\in S^+ $ и $\psi_2\in S_\infty^+. $  Тогда $\psi_1\cdot\psi_2\in S_\infty^+$.
\end{proposition}

   Доказательство  легко провести воспользовавшись  формулой
   \begin{equation*}\label{1l5}
    \frac1{\alpha(\psi_1\cdot\psi_2;t)}= \frac1{\alpha(\psi_1;t)}+ \frac1{\alpha(\psi_2;t)},
   \end{equation*}
справедливость которой вытекает из  цепочки равенств
  $$
   \frac1{\alpha(\psi_1\cdot\psi_2;t)}=-\frac{t(\psi_1(t)\psi_2(t))'}{\psi_1(t)\psi_2(t)}=
      -\frac{t\psi_1'(t)}{\psi_1(t)}-\frac{t\psi_2'(t)}{\psi_2(t)}=\frac1{\alpha(\psi_1;t)}+ \frac1{\alpha(\psi_2;t)}.
  $$

 В качестве примера рассмотрим функцию $\ \psi(t)=t^{-\delta}e^{-\alpha t^r}, \ \ \ \delta\ge0,\ \alpha>0,\ r>0.\ $
Покажем, что она принадлежит множеству $\mathfrak M_\infty^+,$ и для неё  формулы (\ref{8}) и (\ref{9}) являются асимптотическими равенствами при условии выпуклости модуля непрерывности $\omega(t)$.
 Положив $\psi_1(t)=e^{-\alpha t^r}, \ \psi_2(t)= t^{-\delta}$,
несложно убедиться, что $\psi_1\in S_\infty^+$, а $\psi_2\in S^+. $ Поэтому, как следует из утверждения 2, $\psi\in\mathfrak M_\infty^+.$
Для функции $\psi_1(t)$  характеристика $\mu(\psi_1;t)$ имеет вид (\ref{15}).
Выбирая   
$t>\left(\frac{\ln 2}{\alpha{r}}\right)^{{1}/{r}}$, 
получаем
\begin{equation*}\label{16}
  \frac{1}{\mu(\psi_{1};t)}=\sum\limits_{m=1}^{\infty}\frac{\frac{1}{r}\left(\frac{1}{r}-1\right)
  \ldots \left(\frac{1}{r}-m+1\right)}{m!}\left(\frac{\ln 2}{\alpha
  t^{r}}\right)^{m}\leq\frac{A}{ \alpha rt^{r}},
\end{equation*}
где $A$~--- некоторая константа. Из теоремы  \ref{theorem1}, учитывая  предложение 1 и очевидные оценки
\begin{equation*}\label{12_4}
    \frac{\psi(3n)}{\psi(n)}\le  e^{-\alpha n^r( 3^r-1)},
    \ \ \ \  \
    \frac{\psi(n+1)}{\psi(n)}\le \left\{
                                   \begin{array}{ll}
                                     1, & r\in(0,1),\\
                                     e^{- {\alpha r n^{r-1}}}, & r\ge1,
                                   \end{array}
                                 \right.
 \end{equation*}
получаем следующее утверждение.

\begin{corollary}
Пусть $\psi(t)= t^{-\delta}e^{-\alpha t^r} ,\ $ $\delta\ge0,\ \alpha>0, \ r>0,\ $ $\beta\in
\mathbb{R}$ и $\omega(t)$ --- выпуклый модуль
непрерывности.  Тогда при  $n\rightarrow\infty$ выполняются асимптотические равенства
\begin{equation}\label{17}
\left.
  \begin{array}{ll}
    E_n( C_{\beta}^\psi H_{\omega_C})_{C}  \\
    E_n( L_{\beta}^\psi H_{\omega_{L_1}})_{L_1}
  \end{array}
\right\}
  =  \psi(n)\left(\frac{2}{\pi} e_n(\omega)+O(1)\gamma_{n}(\alpha,r)\omega\left(\frac 1n\right)\right),
\end{equation}
где $e_n(\omega)$ определяется равенством (\ref{11}),
а
$\ \gamma_{n}(\alpha,r){=}
   \! \left\{\!\!\!\begin{array}{lc}
 (r\alpha n^{r})^{-1}, &\! r{\in} (0,1),\\
\left(1+ \alpha^{-1}\right) e^{-\alpha} n^{-1}, & r{=}1,\\
  e^{-\alpha r n^{r-1}} n^{-1}, & r{>}1.\\
\end{array}
\right.
$
\end{corollary}

При  $r=1$ и $\delta=0$ либо $r=1$ и $\delta=1$ формулы (\ref{17})     уточняют асимптотические равенства, полученные в  \cite[теорема 1]{Serdyuk_Sokolenko2010_2}, \cite[теорема 1]{Serdyuk_Sokolenko2011_1} и \cite[следствие 1]{Serdyuk_Sokolenko2011_2}. Улучшение достигается за счет более тонкой оценки остаточных членов.

Что же касается случаев $r\in(0,1)$ и $r>1$, то равенства (\ref{17}) являются новыми даже при $\delta=0$. При этом наиболее интересен, на наш взгляд,  случай $r\in(0,1)$, когда применение сумм Фурье либо некоторых иных полиномов, порождаемых классическими линейными методами приближения, для оценок наилучших приближений оказывается неэффективным.

Сопоставляя оценки (\ref{17})   при $\delta=0$ с соответствующими оценками для сумм Фурье  \cite[\S\S5.10-5.16]{Stepanets2002_1}, \cite{Stepanets2001}, видим, что для выпуклых модулей непрерывности $\omega(t)$, выполняются следующие  соотношения:
$$
\frac{ E_n( C_{\beta}^{\alpha,r} H_{\omega_C})_{C}}{{\cal E}  (C_{\beta}^{\alpha,r} H_{\omega_C}; S_{n-1} )_{C}}\sim
\frac{ E_n(L_{\beta}^{\alpha,r} H_{\omega_{L_1}})_{L_1}}{{\cal E}  (L_{\beta}^{\alpha,r} H_{\omega_{L_1}}; S_{n-1} )_{L_1}}\sim
\left\{\begin{array}{l l l}
\displaystyle\frac{\pi}{(1-r)\ln n}, & r\in (0,1),\\
1, & r>1,\\
\end{array}
\right.
$$
$$
\frac{ E_n( C_{\beta}^{\alpha,1} H_{\omega_C})_{C}}{{\cal E}  (C_{\beta}^{\alpha,1} H_{\omega_C}; S_{n-1} )_{C}}\sim \frac{\pi}{2\mathbf{K}\left(e^{-\alpha}\right)}<1,
$$
где $\ \mathbf{K}(q)=\int\limits_{0}^{{\pi}/{2}}(1-q^{2}\sin^{2}u)^{-1/2}du,\ $ $ 0\leq q<1,\ $~--- полный эллиптический интеграл первого рода,
 а запись $A(n){\sim} B(n)$  означает выполнение предельного равенства
$\lim\limits_{n\rightarrow\infty} {A(n)}/{B(n)}{=}1.$

Таким образом, на классах $C_{\beta}^{\alpha,r} H_{\omega_C}$ и $L_{\beta}^{\alpha,r} H_{\omega_{L_1}}$ суммы Фурье реализуют асимптотику наилучших приближений только при $r>1.$  Аналогичный эффект для классов $C_{\beta}^{\alpha,r} B_{L_\infty}$ и $L_{\beta}^{\alpha,r} B_{L_1}$, где $B_{L_p}=\{\varphi\in L_p: \ \|\varphi\|_{L_p}\le 1\},$ был известен ранее \cite{Serdyuk2004}.

\begin{center}
\textbf{Доказательство теоремы \ref{theorem1}}
\end{center}
   Установим сначала оценки сверху в (\ref{8}) и (\ref{9}) и неравенство (\ref{10}). Для этого  будем использовать линейный метод приближения $U^{\ast}_{n-1}$, рассмотренный в работе  \cite{Serdyuk2004}.
Каждой функции $f$ из класса $L_{\beta}^\psi\mathfrak N$ {($C_{\beta}^\psi\mathfrak N$)}  с помощью системы чисел
\begin{equation*}\label{21}
      \begin{array}{l}
              \displaystyle\lambda_{k}^{(n)}=\lambda_{k}^{(n)}(\psi;\beta)=(\psi(k)-\psi(2n-k)-\psi(2n+k))\cos\frac{\beta\pi}{2},  \\
       \displaystyle\nu_{k}^{(n)}=\nu_{k}^{(n)}(\psi;\beta)=(\psi(k)-\psi(2n-k)+\psi(2n+k))\sin\frac{\beta\pi}{2},
     \end{array} \ \ \
         k=1,\ldots,n-1.
\end{equation*}
приведем в соответствие тригонометрический полином
 \begin{equation*}\label{22}
U^{\ast}_{n-1}(f;x)
=\frac{a_0(f)}{2}+\sum\limits_{k=1}^{n-1}\left(\lambda_{k}^{(n)}(a_{k}\cos kx +b_{k}\sin kx)
  + \nu_{k}^{(n)}(a_{k}\sin kx -b_{k}\cos kx)\right),
\end{equation*}
где $a_{k}=a_{k}(f_{\beta}^{\psi})$,
$b_{k}=b_{k}(f_{\beta}^{\psi})$, $k\in \mathbb N$,~--- коэффициенты Фурье функции $f_{\beta}^{\psi}$.

Как следует из леммы 2  работы  \cite{Serdyuk2004}, для произвольной функции $f \in L_{\beta}^{\psi}\mathfrak N$,
$\mathfrak N\subseteq L_1$, почти при всех $x \in \mathbb{R}$ имеет место равенство
\begin{equation}\label{23}
f(x)-U^{\ast}_{n-1}(f;x)=
\frac{1}{\pi}\int\limits_{0}^{2\pi} f^\psi_\beta(x-t)
\bigg(2\cos\left(nt-\frac{\beta\pi}{2}\right) \Psi_{n}(t) -\Psi_{n,\beta}(t)\bigg)dt,
\end{equation}
в котором
\begin{equation*}\label{24}
\Psi_{n}(t)= \frac{\psi(n)}{2} + \sum
\limits_{k=1}^{\infty} \psi(n+k) \cos kt,
\ \ \
    \Psi_{n,\beta}(t) =\sum \limits_{k=n}^{\infty} \psi(2n+k) \cos \left(kt+\frac{\beta\pi}{2}\right).
\end{equation*}
Если же $f\in C_{\beta}^{\psi}\mathfrak N$, то интегральное представление
(\ref{23}) выполняется при всех $x \in \mathbb{R}$.

Введем обозначение $X_{\beta}^{\psi} H_{\omega_X}, $ где  $X$ есть $C$ или $L_p, 1\le p<\infty, $ которое следует понимать следующим образом:
$$
X_{\beta}^{\psi} H_{\omega_X}=\left\{
                                    \begin{array}{ll}
                                      C_{\beta}^{\psi} H_{\omega_C}, & \hbox{если}\ X=C, \\
                                      L_{\beta}^{\psi} H_{\omega_{L_p}}, & \hbox{если}\ X=L_p,\ \  1\le p<\infty.
                                    \end{array}
                                  \right.
$$

Принимая во внимание  представление (\ref{23}),   при $f\in X_{\beta}^{\psi} H_{\omega_X} $
имеем
\begin{equation}\label{25}
E_n(f)_{X}=E_n(f   -U^\ast_{n-1}(f))_{X}=
$$
$$
   =\inf\limits_{T_{n-1}\in {\cal T}_{2n-1}}\left\|  \frac1\pi  \int\limits_{0}^{2\pi}
f^\psi_\beta(\cdot-t) \bigg( 2\cos\left(nt-\frac{\beta\pi}2\right)\Psi_{n}(t)  - \Psi_{n,\beta}(t)   \bigg)dt-T_{n-1}(\cdot)\right\|_{X}.
\end{equation}

Правая часть в (\ref{25}) является $4$-периодической функцией по $\beta.$ Поэтому далее будем считать, что $\beta\in[0,4).$
Учитывая это, положим
\begin{equation}\label{26}
x_k=\frac{(1+\beta)\pi}{2n}+\frac{k\pi}n,\quad
  t_k=x_k-\frac{\pi}{2n}=\frac{\beta\pi}{2n}+\frac{k\pi}n,\quad k\in\mathbb{Z},
\end{equation}
\begin{equation}\label{28}
l_n(t)=x_k  \ \  \hbox{при }\ \ t\in[t_k,t_{k+1}),\ \ \ k\in\mathbb{Z},
\end{equation}
и рассмотрим кусочно-постоянную $2\pi$-периодическую функцию $\overline{\Psi}_{n}(t)$ вида
\begin{equation}\label{29}
    \overline{\Psi}_{n}(t)=\Psi_{n}(l_n(t)).
\end{equation}

Из равенства (\ref{25}),
применяя неравенство треугольника для нормы, обобщенное неравенство Минковского \cite[c. 395]{Kornejchuk1987}
\begin{equation}\label{31}
 \left\|\int\limits_0^{2\pi} F(\cdot,t) dt\right\|_{L_p}\le\int\limits_0^{2\pi}\left\|F(\cdot,t)\right\|_{L_p} dt , \ \ \ 1\le p\le \infty,
\end{equation}
и  учитывая тот факт, что свертка любой функции из $L_1$ с тригонометрическим полиномом из ${\cal T}_{2n-1}$ также является полиномом из ${\cal T}_{2n-1}$, при $f\in X_{\beta}^{\psi} H_{\omega_X} $ 
получаем
$$
    E_n(f)_{X}
\le
 \frac2\pi \left\| \int\limits_{0}^{2\pi}
f^\psi_\beta(\cdot-t)  \cos\left(nt-\frac{\beta\pi}2\right) \overline{\Psi}_{n}(t) dt\right\|_{X}+
$$
$$
{+} \inf\limits_{T_{n{-}1}\in {\cal T}_{2n{-}1}}\!\left\|\int\limits_{0}^{2\pi}
\left(f^\psi_\beta(\cdot-t){-}T_{n{-}1}(\cdot{-}t)\right)\!\!   \left(\cos\!\left(\!nt{-}\frac{\beta\pi}2\!\right)\! \left({\Psi}_{n}(t){-}\overline{\Psi}_{n}(t)\right){-}{\Psi}_{n,\beta}(t)\!\right)\!dt\right\|_{X}{\le}
$$
\begin{equation}\label{32}
{\le} \frac2\pi \left\| \int\limits_{0}^{2\pi}
f^\psi_\beta(\cdot{-}t)  \cos\left(nt{-}\frac{\beta\pi}2\right) \overline{\Psi}_{n}(t) dt\right\|_{X}+
E_{n}(f^\psi_\beta)_{X}
 \int\limits_{0}^{2\pi}\left(|{\Psi}_{n}(t){-}\overline{\Psi}_{n}(t)|{+}
| {\Psi}_{n,\beta}(t) |\right)dt.
\end{equation}

Установим  оценку сверху для последнего интеграла в (\ref{32}).
 Учитывая (\ref{26})--(\ref{29}), имеем
\begin{equation}\label{35}
\int
    \limits_{0}^{2\pi}\left|{\Psi}_{n}(t){-}\overline{\Psi}_{n}(t)\right|dt=\sum_{k=0}^{2n-1}\int\limits_{t_k}^{t_{k+1}}\left|{\Psi}_{n}(t){-}
\overline{\Psi}_{n}(t)\right|dt\le
\frac\pi n\sum_{k=0}^{2n-1}\mathop{\rm V}_{t_k}^{t_{k+1}}(\Psi_n)=\frac{\pi}n\mathop{\rm V}_0^{2\pi} (\Psi_n).
\end{equation}

Оценим вариацию функции $\Psi_{n}(t)$ на $[-\pi,\pi]$. Если
$\psi\in \mathfrak M_\infty^+$, то, как следует из \cite[c. 288--289]{Stepanets2002_1}, ядро $\Psi_{n}(t)$ можно дифференцировать сколько угодно раз и ряды, получаемые почленным дифференцированием, будут равномерно сходиться на $[0, 2\pi]$. Поэтому
\begin{equation}\label{91}
\mathop{\mathrm{V}}\limits^{2\pi}_{0}(\Psi_{n})=\int\limits_{0}^{2\pi}
|\Psi_{n}'(t)|dt=\int\limits_{0}^{2\pi}\left\vert\sum\limits_{k=1}^{\infty}
k\psi(n+k)\sin kt\right\vert dt.
\end{equation}
Для оценки интеграла в правой части (\ref{91})
воспользуемся следующим   результатом Теляковского \cite{Telyakovskii1964}.

\textsc{Теорема  Теляковского.} \textit{Пусть
последовательность $\{\alpha_{k}\}_{k=0}^{\infty}$
удовлетворяет условиям
\begin{equation}\label{92}
  \lim\limits_{k\rightarrow\infty}\alpha_{k}=0, \ \ \ \ \ \sum\limits_{k=0}^{\infty}
  (k+1)\left\vert\Delta^{2}\alpha_{k}\right\vert<\infty, \ \ \ \ \
  \sum\limits_{k=1}^{\infty}\frac{|\alpha_{k}|}{k}<\infty.
\end{equation}
Тогда тригонометрический ряд $\ \sum\limits_{k=0}^{\infty}\alpha_{k}\sin
kt\ $ является рядом Фурье некоторой функции $g\in L_1$ и
\begin{equation}\label{95}
  \|g\|_{L_1}\leq A\left(\sum\limits_{k=0}^{\infty}(k+1)\left|
  \Delta^{2}\alpha_{k}\right|+
  \sum\limits_{k=1}^{\infty}\frac{|\alpha_{k}|}{k}\right),
\end{equation}
 где $A$~---   абсолютная константа.}

Положим
$\alpha_{k}=k\psi\left(n+k\right)$ и докажем выполнение условий
(\ref{92}) теоремы Теляковского в случае, когда $\psi\in\mathfrak M_\infty^+$. Используя  неравенство 3.12.24 из \cite{Stepanets2002_1}
\begin{equation}\label{80}
  \psi(t)\leq -2\psi'(t)\left(\eta(t)-t)\right), \ \ t\geq 1, \ \ \psi\in
  \mathfrak M_\infty^+,
\end{equation}
имеем
$$
 \left(t \psi(n+t)\right)'\le t |\psi'(n+t)|\left(\frac{2(n+t)}{t\mu(n+t)}-1\right)\le
t \left\vert\psi'(n+t)\right\vert\left(\frac{2n+2}{\mu(n+t)}-1\right).
$$
Величина $(2n+2)/\mu(n+t)$ монотонно стремится к 0
при $t\rightarrow\infty,$ а значит
\begin{equation}\label{97'}
 \mathop{\lim}\limits_{k\rightarrow\infty}\alpha_{k}=
 \mathop{\lim}\limits_{k\rightarrow\infty}k\psi(n+k)=0.
\end{equation}

Докажем, сходимость ряда
 $
\sum\limits_{k=0}^{\infty}(k+1)|\Delta^{2}\alpha_{k}|=\sum\limits_{k=0}^{\infty}(k+1)|\Delta^{2}
\{k\psi(n+k)\}|.
 $
Поскольку
 $
\Delta^{2}\{k\psi(n+k)\}=
k\Delta^{2}\psi(n+k)-2\Delta\psi(n+k+1),
 $
то, вследствие выпуклости вниз  последовательности $\psi(k)$, имеем
 $$
\sum\limits_{k=0}^{\infty}(k+1)\left\vert\Delta^{2}\{k\psi(n+k)\}\right\vert
\leq
\sum\limits_{k=0}^{\infty}(k+1)k\Delta^{2}\psi(n+k)+\sum\limits_{k=0}^{\infty}
2(k+1)\Delta\psi(n+k+1)=
 $$
\begin{equation}\label{98}
=\sum\limits_{k=1}^{\infty}k^{2}\Delta^{2}\psi(n+k)+\psi(n+1)+2\sum\limits_{k=1}^{\infty}\psi(n+k).
\end{equation}

Для оценки первой суммы в правой части (\ref{98})
воспользуемся преобразованием Абеля
\begin{equation}\label{99}
  \sum\limits_{k=1}^{N}u_{k}v_{k}=\sum\limits_{k=1}^{N-1}(u_{k}-u_{k+1})V_{k}+u_{N}V_{N},
\end{equation}
где   $\ V_{k}=v_{1}+v_{2}+\ldots+v_{k}.\ $   Дважды применяя преобразование Абеля (\ref{99}), получаем
 \begin{equation}\label{100}
\sum\limits_{k=1}^N k^{2}\Delta^{2}\psi(n{+}k)
  =\psi(n+1)+2\sum\limits_{k=1}^{N-1}\psi(n+k+1)-(2N-1)\psi(n+N+1)-N^{2}
  \Delta\psi(n+N+1).
\end{equation}
Из равенства (\ref{100}) вытекает
\begin{equation}\label{101}
  \sum\limits_{k=1}^{N}k^{2}\Delta^{2}\psi(n+k)<\psi(n+1)+2\sum\limits_{k=2}^{N}\psi(n+k),
  \ \ N\in \mathbb{N}.
\end{equation}
Переходя к пределам при $N\rightarrow\infty$ в обеих частях
неравенства (\ref{101}), получаем
\begin{equation}\label{102}
  \sum\limits_{k=1}^{\infty}k^{2}\Delta^{2}\psi(n+k)\leq\psi(n+1)+2\sum\limits_{k=2}^{\infty}\psi(n+k).
\end{equation}
Сопоставление формул (\ref{98}) и (\ref{102}) позволяет записать соотношение
\begin{equation}\label{103}
\sum\limits_{k=0}^{\infty}(k+1)|\Delta^{2}\alpha_{k}|=\sum\limits_{k=0}^{\infty}
(k+1)\left\vert\Delta^{2}\{k\psi(n+k)\}\right\vert\leq
 4\sum\limits_{k=1}^{\infty}\psi(n+k)<\infty.
  \end{equation}
В силу определения последовательности $\alpha_{k}$ и включения
\mbox{$\psi \in \mathfrak M_\infty^+$}
\begin{equation}\label{104}
  \sum\limits_{k=1}^{\infty}\frac{|\alpha_{k}|}{k}=\sum\limits_{k=1}^{\infty}\psi(n+1)<\infty.
\end{equation}
Как видно из (\ref{97'}), (\ref{103}) и (\ref{104}), при $\alpha_{k}=k\psi(n+k)$ все условия
теоремы  Теляковского выполняются. Поэтому на основании  (\ref{95}) с учетом (\ref{103}) и (\ref{104}) имеем
\begin{equation}\label{105}
\int\limits_{0}^{2\pi}\left\vert\sum\limits_{k=1}^{\infty}k\psi(n+k)\sin
kt \right\vert dt\le A\sum\limits_{k=1}^{\infty}\psi(n+k).
\end{equation}
 Из (\ref{91}) и
(\ref{105}) следует оценка
\begin{equation}\label{106}
  \mathop{\mathrm{V}}\limits_{0}^{2\pi}\left(\Psi_{n}\right)=
  O(1)\sum\limits_{k=1}^{\infty}\psi(n+k).
\end{equation}

Поскольку $\psi\in \mathfrak M_\infty^+$, то из формулы (19) работы \cite{Serdyuk2004_2} вытекает, что для всех $m\geq 1$  таких, что $\mu(m)>2$,
выполняется неравенство
\begin{equation}\label{107}
  \int\limits_{m}^{\infty}\psi(v)dv\leq
  \frac{2}{1-\frac{2}{\mu(m)}}\psi(m)(\eta(m)-m), \ \ m\geq 1.
\end{equation}
Применяя (\ref{107}) при $m=n+1$, имеем
 \begin{equation}\label{108}
\sum\limits_{k=1}^{\infty}\psi(n+k)\leq
\psi(n+1)+\int\limits_{n+1}^{\infty}\psi(v)dv
\leq
\psi(n+1)\Big(1+\frac{2}{1-\frac{2}{\mu(n+1)}}(\eta(n+1)-(n+1))\Big).
\end{equation}

Сопоставляя соотношения (\ref{35}), (\ref{106}) и (\ref{108}), получаем
\begin{equation}\label{34}
     \int    \limits_{0}^{2\pi}\left|{\Psi}_{n}(t){-}\overline{\Psi}_{n}(t)\right|dt
     =O(1)\psi(n+1)\left(\frac1{n}+\frac{1}{\mu(n)-2}\right).
\end{equation}

Для окончательной оценки последнего интеграла в (\ref{32}) покажем, что
\begin{equation}\label{83}
  \int\limits_{0}^{2\pi}| {\Psi}_{n,\beta}(t) |dt=O(1)\psi(3n)\left(1+\ln^{+}(\eta(n)-n)\right).
\end{equation}
Для сумм $K_n(t)=\sum\limits_{k=n}^\infty\alpha_k\cos(kt+\delta), \ n\in\mathbb N,\ \delta\in\mathbb R,\ \alpha_k\downarrow0,\ \Delta^2\alpha_k\ge0, \ $ хорошо известна оценка \cite[c. 135]{Stechkin1980}
\begin{equation}\label{83'}
    \|K_{n}\|_{L_1}=Q_{n}+O(\alpha_{n}),\ \ \
Q_{n}=
  \left\{
\begin{array}{l l}
\displaystyle\frac{4}{\pi}\sum\limits_{k=n}^{2n-1}\frac{\alpha_{k}}{k+1-n}+2|\sin\delta|\sum\limits_{k=2n}^{\infty}
\frac{\alpha_{k}}{k}, & \delta\neq\pi l,  l\in \mathbb{Z},\\
\displaystyle\frac{4}{\pi}\sum\limits_{k=n}^{2n-1}\frac{\alpha_{k}}{k+1-n},& \delta=\pi
l,  l\in \mathbb{Z},
\end{array}\right.
\end{equation}
(равенство (\ref{83'}) является  следствием результатов работы Теляковского \cite{Telyakovskii1971}). При  $\alpha_{k}=\psi(2n+k)$ и $\delta=\beta\pi/2$ из (\ref{83'})
получаем
\begin{equation}\label{77}
\frac12\int\limits_{-\pi}^{\pi}\left|\Psi_{n,\beta}(t)\right|
  dt
<   \sum\limits_{k=n}^{\infty}\frac{\psi(2n+k)}{k+1-n}\le
  \psi(3n)+  \int\limits_{1}^{T(n)}
 \frac{\psi(t + 3n - 1)}{t} dt
  +  \int\limits_{T(n)}^{\infty} \frac{\psi(t + 3n - 1)}{t} dt,
\end{equation}
где $ T(n)=\max\{1; \eta(n)-n\},  \ \ n \in
\mathbb{N}. $

Оценим каждый из интегралов в (\ref {77}).
Легко видеть, что
\begin{equation}\label{79}
\int\limits_{1}^{T(n)}\frac{\psi(t+3n-1)}{t}\,dt\leq
\psi(3n)\int\limits_{1}^{T(n)}\frac{dt}{t}= 
\psi(3n)\ln^{+}\left(\eta(n)-n\right).
\end{equation}
Используя   факт монотонного убывания функции $\xi(t)=(\eta(a+t)-(a+t))/{t}$ на $[1,\infty)$ при \mbox{$a>0$}  и применяя неравенство (\ref {80}),
 получаем
 \begin{equation}\label{81}
\int\limits_{T(n)}^{\infty}\frac{\psi(t+3n-1)}{t}dt \leq
  -2\int\limits_{T(n)}^{\infty}\psi'(t+3n-1)\frac{\eta(t+3n-1)-(t+3n-1)}{t}dt
\le
 $$
 $$
 \le  -2
\frac{\eta\left(T(n)+3n-1\right)-\left(T(n)+3n-1\right)}{T(n)+3n-1}\cdot\frac{T(n)+3n-1}{T(n)}\int\limits_{T(n)}^{\infty}\psi'(t+3n-1)dt\le
 $$
 $$
\le\frac{2(1+ 3\mu(n))}{\mu(3n)}\psi(T(n)+3n-1)
\le 2\left(3+\frac1{\mu(n)}\right)\psi(3n).
\end{equation}
Объединяя  (\ref{77})--(\ref{81})   убеждаемся в истинности (\ref{83}).

Из   (\ref{32}),   (\ref{34}), (\ref{83})
 и     неравенства Джексона $\
      E_n(\varphi)_X\le \displaystyle\frac 32\omega\left(\varphi,\frac \pi n\right)_X,\   \varphi\in X$ (см., например, \cite[c. 258]{Kornejchuk1987}),
   получаем
\begin{equation}\label{30}
     E_n(f)_{X}\le\frac2\pi\left\|  \int\limits_{0}^{2\pi}
f^\psi_\beta(\cdot-t)  \cos\left(nt-\frac{\beta\pi}2\right) \overline{\Psi}_{n}(t) dt\right\|_{X}+O(1)\gamma_n(\psi)\omega\left(\frac 1n\right),
 \end{equation}
где  $\gamma_n(\psi)$  определяется равенством (\ref{12}).

Рассматривая точные верхние грани по $f\in X_{\beta}^{\psi} H_{\omega_X} $
в обеих частях неравенства (\ref{30}), получаем
 \begin{equation}\label{37}
     E_n(X_{\beta}^{\psi} H_{\omega_X})_{X}\le I_n(\psi,\beta,\omega)_{X}+O(1)\gamma_n(\psi)\omega\left(\frac 1n\right),
 \end{equation}
где
$\  I_n(\psi,\beta,\omega)_{X}=\sup\limits_{\varphi\in H_{\omega_X}}I_n(\varphi)_{X}, \ \  I_n(\varphi)_{X}= \displaystyle\frac2\pi\left\|  \int\limits_{0}^{2\pi}
\varphi(\cdot-t)  \cos\left(nt-\frac{\beta\pi}2\right) \overline{\Psi}_{n}(t) dt\right\|_{X}.$

Для оценки величины $I_n(\psi,\beta,\omega)_{X} $  нам будет необходимо такое утверждение.

\begin{lemma}\label{lemma_2}  Пусть $X$ есть $C$ или
$L_p,\ 1\le p<\infty,\ \psi\in\mathfrak M_\infty^+,\ $  $\ \beta
\in\mathbb R\ $  и $\ \omega (t)$ --- произвольный модуль непрерывности. Тогда при
$\ \mu(n)>2$  имеет место неравенство
\begin{equation}\label{39}
I_n(\psi,\beta,\omega)_{X} \le \frac{2
}{\pi}\psi(n)e_n(\omega)
+O(1)   \psi(n+1)\left(\frac1n+\frac1{\mu(n)-2}\right)\omega\left(\frac 1n\right),
\end{equation}
в котором $e_n(\omega)$ определяется равенством (\ref{11}).
\end{lemma}

{\scshape Доказательство леммы \ref{lemma_2}.} Учитывая обозначения (\ref{26})--(\ref{29}), а также положительность функции $\Psi_n,$ для произвольной   $\varphi\in H_{\omega_X}$
имеем
\begin{equation}\label{40}
    I_n(\varphi)_{X}
= \frac2\pi \left\|\sum\limits_{k=0}^{2n-1}\Psi_n(x_k)
\int\limits_{t_k}^{t_{k+1}}\varphi(\cdot-t) \cos\left(n
t-\frac{\beta\pi}2\right) dt \right\|_{X}
\le\frac2\pi\sum\limits_{k=0}^{2n-1}\Psi_n(x_k) e_k(\varphi,n)_X,
\end{equation}
где $\  e_k(\varphi,n)_X=\left\|\int\limits_{t_k}^{t_{k+1}}\varphi(\cdot-t)\cos\left(n
t-\frac{\beta\pi}2\right)dt\right\|_X.\
$
Применяя неравенство Минковского (\ref{31}),
получаем
$$
    e_k(\varphi,n)_X
   =  \left\|\int\limits_{t_k}^{x_{k}}(\varphi(\cdot-t)-\varphi(\cdot-2x_k+t))\cos\left(n
t-\frac{\beta\pi}2\right)dt\right\|_X\le
$$
\begin{equation}\label{42}
{\le}
  \int\limits_{t_k}^{x_{k}}\omega(2(x_k-t))\left|\cos\left(n
t-\frac{\beta\pi}2\right)\right|dt=
\int\limits_0^{\pi/2n}\omega(2t)\sin n t
dt=\frac1n\int\limits_0^{\pi/2}\omega\left(\frac{2t}n\right)\sin t
dt.
\end{equation}

Объединяя соотношения (\ref{40}) и (\ref{42}),   можем записать
\begin{equation}\label{43}
      I_n(\psi,\beta,\omega)_X\leq \frac2{\pi n}\int\limits_0^{\pi/2}\omega\left(\frac{2t}n\right)\sin t
dt\sum\limits_{k=0}^{2n-1} {\Psi_n(x_k)} .
\end{equation}

Покажем, что
\begin{equation}\label{44}
    \frac1n\sum\limits_{k=0}^{2n-1} {\Psi_n(x_k)}  = \psi(n)
    + O(1)\psi(n+1)\left(\frac1n+\frac1{\mu(n)-2}\right).
\end{equation}

Действительно,  в силу (\ref{29}) и (\ref{34})
$$
  \frac1n  \sum\limits_{k=0}^{2n{-}1}{\Psi_n(x_k)} {=}\sum\limits_{k=0}^{2n{-}1}\frac{1}{\pi}\int\limits_{t_k}^{t_{k{+}1}} \overline{\Psi}_n(t)dt{=}
 \frac{1}{\pi}\sum\limits_{k=0}^{2n{-}1}\int\limits_{t_k}^{t_{k{+}1}}\Psi_n(t)dt
    +\frac{1}{\pi}\sum\limits_{k=0}^{2n{-}1}\int\limits_{t_k}^{t_{k{+}1}}(\overline{\Psi}_n(t){-}\Psi_n(t))dt{=}
$$
$$
= \frac1\pi\int\limits_0^{2\pi}\Psi_n(t)dt +O(1)\int\limits_{0}^{2\pi}|\overline{\Psi}_n(t)-\Psi_n(t)|dt=
\psi(n)
    +  O(1)\psi(n+1)\left(\frac1n+\frac1{\mu(n)-2}\right).
$$

Объединяя (\ref{43}) и (\ref{44}),  получаем неравенство (\ref{39}). Лемма 1 доказана.

Продолжим доказательство теоремы 1. Из сопоставления соотношений  (\ref{37})  и (\ref{39}), следуют оценки сверху в (\ref{8})--(\ref{10}).

Докажем  оценку снизу в (\ref{8}). Пусть   $\omega(t)$
--- выпуклый модуль непрерывности. Рассмотрим функцию
$f^\ast(x)$,
связанную равенством (\ref{1}) с функцией
\begin{equation*}\label{45}
\varphi^\ast(t)= (-1)^{k}\varphi_k(t),\ \ t\in [t_k,t_{k+1}], \
\ k\in\mathbb{Z},
\end{equation*}
где
\begin{equation*}\label{46}
    \varphi_k(t)=\left\{%
\begin{array}{rl}
     \omega\big(2(x_k-t)\big)/2, & t\in [t_k,x_k], \\
    - \omega\big(2(t-x_k)\big)/2, & t\in [x_k,t_{k+1}], \\
\end{array}%
\right.\ \ \ k\in\mathbb{Z}.
\end{equation*}
По построению   $\varphi^\ast(t)$   является непрерывной  $2\pi$-периоди\-ческой функцией,  $\varphi^\ast\bot 1$ и, в силу выпуклости   $\omega(t)$, удовлетворяет неравенству
\begin{equation}\label{47}
    |\varphi^\ast(t_1)-\varphi^\ast(t_2)|\leq\omega(|t_1-t_2|) \ \ \
    \forall t_1,t_2\in \mathbb{R}.
\end{equation}
Значит,  $f^\ast\in C^{\psi}_\beta H_{\omega_C}.$  Найдем оценку
величины  наилучшего равномерного приближения   функции $f^\ast(x) $.
 С учетом  (\ref{23}) и (\ref{29}),
  можем записать
\begin{equation}\label{48'}
E_n(f^\ast)_C=E_n(f^\ast -U^\ast_{n-1}(f^\ast))_C =E_n(\Phi^\ast_n)_C+O(1)\|r_n\|_C,
\end{equation}
где
\begin{equation}\label{48}
\Phi^\ast_n(x)=\Phi^\ast_n(\varphi^\ast;x)=\frac2\pi  \int\limits_{0}^{2\pi}
\varphi^\ast(x-t)  \cos\left(nt-\frac{\beta\pi}2\right) \overline{\Psi}_{n}(t) dt,
\end{equation}
\begin{equation}\label{49}
 r_n(x) = r_n(\varphi^\ast;x)=  \int\limits_{0}^{2\pi}\varphi^\ast(\cdot-t)
   \bigg(\cos\left(nt-\frac{\beta\pi}2\right) \left({\Psi}_{n}(t)-\overline{\Psi}_{n}(t)\right)+{\Psi}_{n,\beta}(t)\bigg) dt  .
\end{equation}

Применяя неравенство Минковского (\ref{31}) к правой части (\ref{49}), а затем учитывая оценки (\ref{34}), (\ref{83}) и определение функции $\varphi^\ast(t)$, получаем
\begin{equation}\label{51}
    \|r_n\|_C=O(1)\gamma_n(\psi)\omega\left(\frac 1n\right).
\end{equation}

Из построения  $\varphi^\ast(t)$ следует, что в $2n$ точках $\zeta_k=  {k\pi}/ n, \ k=0,1,\ldots,2n-1,$ функция $\Phi^\ast_n(x)$ принимает значения с чередующимися
знаками:
$$\
{\rm sign}\,\Phi^\ast_n(\zeta_k)=-{\rm sign}\,\Phi^\ast_n(\zeta_{k+1}),
\ \ \ k=0,1,\ldots, 2n-2.
$$

Кроме того, при $k=0,1,\ldots, 2n-1$
\begin{equation*}\label{50}
    |\Phi^\ast_n(\zeta_k)|=\frac2{\pi n}\left|\sum\limits_{k=0}^{2n-1} {\Psi_n(x_k)}   \int\limits_{t_k}^{t_{k+1}}\varphi^\ast(\zeta_k-t)\cos \left(nt-\frac{\beta\pi}2\right)dt\right|=
$$
$$
{=}
\frac2{\pi n}\sum\limits_{k=0}^{2n-1} {\Psi_n(x_k)}   \int\limits_{t_k}^{x_k}\omega(2(x_k-t))\left|\cos \left(nt{-}\frac{\beta\pi}2\right)\right|dt{=}
\frac2{\pi n}\int\limits_0^{\pi/2}\omega\left(\frac{2t}n\right)\sin t
dt\sum\limits_{k=0}^{2n-1} {\Psi_n(x_k)}.
\end{equation*}
Используя (\ref{44}) и применяя  теорему Чебышева об альтернансе (см., например, \cite[теорема 2.1.2]{Kornejchuk1987}),  находим
\begin{equation}\label{51'}
E_n(\Phi^\ast_n)_C=|\Phi^\ast_n(\zeta_k)|=    \frac{2}{\pi} \psi(n)e_n(\omega)+O(1)\psi(n+1)\left(\frac1n+\frac1{\mu(n)-2}\right)\omega\left(\frac 1n\right).
\end{equation}
Учитывая (\ref{48'}), (\ref{51})  и (\ref{51'}), для выпуклых модулей непрерывности $\omega(t)$ имеем
\begin{equation}\label{52}
 E_n(C^{\psi}_\beta H_{\omega_C})_C\ge E_n(f^\ast)_C\ge  \frac{2}{\pi} \psi(n)e_n(\omega)+O(1)\gamma_{n}(\psi)\omega\left(\frac 1n\right),
\end{equation}
где $e_n(\omega)$ определяется равенством (\ref{11}), а $\gamma_{n}(\psi)$ --- равенством (\ref{12}).

Если же $\omega(t)$ --- произвольный  модуль непрерывности, то  для функции $\varphi^\ast(t)$ условие (\ref{47}) может не
выполняться, однако как известно (см., например, \cite[c. 202]{Stepanets2002_1})     функция
$\varphi_\ast(t)=2\varphi^\ast(t)/3$   принадлежит к классу
$H_{\omega_C}.$ Следовательно $f_\ast(x)= 2 f^\ast (x)/3 $   принадлежит к  $C^{\psi}_\beta
H_{\omega_C}  $ и
\begin{equation}\label{53}
   E_n(C^{\psi}_\beta H_{\omega_C})_C\ge E_n(f_\ast)_C=\frac23E_n(f^\ast)_C\ge
\frac{4}{3\pi} \psi(n)e_n(\omega)+O(1)\gamma_{n}(\psi)\omega\left(\frac 1n\right).
\end{equation}

Объединяя (\ref{37}),  (\ref{39}), (\ref{52}) и (\ref{53}), получаем соотношение
(\ref{8}) для произвольных модулей непрерывности $\omega(t)$.

Установим необходимую оценку снизу в (\ref{9}).
Пусть  $\omega(t)$
--- выпуклый  модуль непрерывности.
Положим
\begin{equation*}\label{54}
    \varphi_1(t)=\left\{%
\displaystyle\begin{array}{rl}
    \displaystyle \omega(2t)/4, & t\in [0, \pi/{2n}), \\
    -\displaystyle \omega(-2t)/4, & t\in (-\pi/{2n},0], \\
     0, &   \pi/{2n}\le |t|\le \pi, \\
\end{array}%
\right.
\end{equation*}
и через $\varphi_2(t)$ обозначим $2\pi$-периодическое продолжение функции $\varphi_1(t)$. Далее, рассмотрим   функцию
\begin{equation}\label{55}
    \varphi^\ast(t)=\varphi^\ast_\omega(t)=\varphi_2'(t)-\frac1{4\pi}\omega\left(\frac\pi n\right).
\end{equation}
Как известно (см., например,  \cite[c.~258]{Stepanets2002_1}), $\varphi^\ast\bot1$ и   $\varphi^\ast\in H_{\omega_{L_1}}$.
Для функции $f^\ast\in L^\psi_\beta H_{\omega_{L_1}},$ связанной равенством
(\ref{1}) с функцией $\varphi^\ast(t)$, учитывая  (\ref{23}) и (\ref{29}), можем записать
\begin{equation*}\label{56}
    E_n(f^\ast)_{L_1}=E_n(f^\ast-U_{n-1}^{\ast})_{L_1}=E_n(\Phi^\ast_n)_{L_1}+ O(1)\|r_n\|_{L_1} ,
\end{equation*}
где $\Phi^\ast_n(x)$ и $r_n(x)$  определены равенствами (\ref{48}) и (\ref{49}), соответсвенно.
Для $\|r_n\|_{L_1}$ легко получить аналог  оценки (\ref{51}) вида
$
    \|r_n\|_{L_1}=O(1)\gamma_n(\psi)\omega (  1/n).
$

Функция $ \Phi^\ast_n(x)$ в точках $x_i=\frac{(1+\beta)\pi}{2n}+\frac{i\pi}n,\ \ i\in \mathbb{Z},$ равна нулю
\begin{equation*}\label{60}
 \Phi^\ast_n(x_i){=}\frac2\pi
\int\limits_0^{2\pi}\varphi^\ast(x_i-t)\cos\left(n
t-\frac{\beta\pi}2\right)\overline{\Psi}_n(t)dt=
$$
$$
=\frac{2\overline{\Psi}_n(x_i)}\pi\int\limits_{t_i}^{t_{i+1}}\varphi^\ast(x_i-t)\cos\left(n
t-\frac{\beta\pi}2\right)dt=
 \frac{(-1)^i 2\overline{\Psi}_n(x_i)}\pi\int\limits_{-\pi/2n}^{\pi/2n}\varphi^\ast(t)\sin nt dt=0
\end{equation*}
и других нулей не имеет. При этом
\begin{equation}\label{61}
 {\rm sign}\, \Phi^\ast_n(x)=(-1)^i,\ \ x\in(x_i,x_{i+1}),\ \ i\in\mathbb{Z}.
\end{equation}
Поскольку для функции $\Phi^\ast_n(x)$ выполняется (\ref{61}),
то, как следует из
\cite[теорема 1.4.5]{Kornejchuk1987},\ полином $T_{n-1}^\ast(x)\equiv0$ является полиномом наилучшего приближения в среднем, т.е. $
E_n(\Phi^\ast_n)_{L_1}=\|\Phi^\ast_n\|_{L_1}$.
Значит
\begin{equation}\label{62}
 E_n(L^\psi_\beta H_{\omega_{L_1}})_{L_1}\ge E_n(f^\ast)_{L_1}
=\|\Phi^\ast_n\|_{L_1}+O(1)\gamma_n(\psi)\omega\left(\frac 1n\right).
\end{equation}

Покажем с помощью
стандартных в таких случаях соображений, что
\begin{equation}\label{63}
     \|\Phi^\ast_n\|_{L_1}=  \frac{2}{\pi} \psi(n)e_n(\omega)   +O(1)\psi(n+1)\left(\frac1n+\frac1{\mu(n)-2}\right)\omega\left(\frac 1n\right).
\end{equation}

В силу  (\ref{61})
$$
    \|\Phi^\ast_n\|_{L_1}=\int\limits_{-\pi}^\pi|\Phi^\ast_n(x)|dx=\frac2\pi\sum\limits_{i=0}^{2n-1}
(-1)^i\int\limits_{x_i}^{x_{i+1}}
\int\limits_{t_0}^{t_{2n}}\varphi^\ast(x-t)\cos\left(n
t-\frac{\beta\pi}2\right)\overline{\Psi}_n(t)dtdx=
$$
$$
=
\frac2\pi\sum\limits_{i=0}^{2n-1}(-1)^i
\int\limits_{t_0}^{t_{2n}}\cos\left(n
t-\frac{\beta\pi}2\right)\overline{\Psi}_n(t)   \int\limits_{x_i}^{x_{i+1}}
\varphi^\ast(x-t)dx dt.
$$

Поскольку
$
\int\limits_{x_i}^{x_{i+1}}
\varphi^\ast(x-t)dx=\varphi_2(x_{i+1}-t)-\varphi_2(x_i-t)-\frac1{4n}\omega\left(\frac\pi n\right),
$
то
\begin{equation}\label{65}
    \|\Phi^\ast_n\|_{L_1}{=}\frac2\pi
\sum\limits_{i=0}^{2n-1}({-}1)^i
\int\limits_{t_0}^{t_{2n}}\left(\varphi_2(x_{i+1}{-}t){-}\varphi_2(x_i{-}t){-}\frac1{4n}\omega\left(\frac\pi n\right)\right)\cos\left(n
t-\frac{\beta\pi}2\right)\overline{\Psi}_n(t)dt   =
$$
$$
=\frac2\pi\sum\limits_{k=0}^{2n-1} \Psi_n(x_k)
\int\limits_{t_k}^{t_{k+1}}\sum\limits_{i=0}^{2n-1}(-1)^i(\varphi_2(x_{i+1}-t)-\varphi_2(x_i-t))\cos\left(n
t-\frac{\beta\pi}2\right)dt.
\end{equation}
В силу того, что $
\int\limits_{t_k}^{t_{k+1}}\varphi_2(x_i-t)\cos\left(n
t-\frac{\beta\pi}2\right)dt=0, \ \ i\neq k,
$
имеем
$$
\int\limits_{t_k}^{t_{k+1}}\sum\limits_{i=0}^{2n-1}(-1)^i\varphi_2(x_{i+1}-t)\cos\left(n
t-\frac{\beta\pi}2\right)dt=(-1)^{k-1}\int\limits_{t_k}^{t_{k+1}}\varphi_2(x_k-t)\cos\left(n
t-\frac{\beta\pi}2\right)dt,
$$
$$
-\int\limits_{t_k}^{t_{k+1}}\sum\limits_{i=0}^{2n-1}(-1)^i\varphi_2(x_i-t)\cos\left(n
t-\frac{\beta\pi}2\right)dt=(-1)^{k-1}\int\limits_{t_k}^{t_{k+1}}\varphi_2(x_k-t)\cos\left(n
t-\frac{\beta\pi}2\right)dt.
$$

Подставляя эти выражения в (\ref{65}) и учитывая (\ref{44}), находим
$$
    \|\Phi^\ast_n\|_{L_1}=\frac4\pi\sum\limits_{k=0}^{2n-1} (-1)^{k-1} \Psi_n(x_k) \int\limits_{t_k}^{t_{k+1}}\varphi_2(x_k-t)\cos\left(n
t-\frac{\beta\pi}2\right)dt=
$$
$$
=\frac4\pi\int\limits_{-\pi/2n}^{\pi/2n}\varphi_2(t)\sin nt
dt\sum\limits_{k=0}^{2n-1}{\Psi_n(x_k)}
 {=}\frac{2}{\pi} \psi(n)e_n(\omega)   {+}O(1)\psi(n{+}1)\!\left(\frac1n{+}\frac1{\mu(n){-}2}\right)\!\omega\left(\frac 1n\right)\!.
$$
Равенство (\ref{63}) доказано.

Таким образом, объединяя соотношения   (\ref{62}) и (\ref{63}), получаем, что  для выпуклых модулей непрерывности $\omega(t)$
\begin{equation}\label{70}
 E_n(L^\psi_\beta H_{\omega_{L_1}})_{L_1}\ge E_n(f^\ast)_{L_1}\ge
\frac{2}{\pi} \psi(n)e_n(\omega)   +O(1)\gamma_n(\psi)\omega\left(\frac 1n\right).
\end{equation}

Пусть теперь $\omega(t)$ --- произвольный  модуль непрерывности. Для
построения функции $\varphi^\ast(t)$  воспользуемся результатом
 Стечкина (см., например,  \cite[лемма 3.1.1]{Stepanets2002_1}),
согласно которому для произвольного модуля непрерывности $\omega(t)$ существует
выпуклый модуль непрерывности $\omega^\ast(t)$ такой, что
$\
\omega(t)\le \omega^\ast(t)<2\omega(t)\ \ \ \forall t>0.
$

Поскольку $\bar\omega(t)= \omega^\ast(t)/2$ --- выпуклая функция,
то построив по приведенной выше схеме функцию
$\varphi^\ast(t)=\varphi^\ast_{\bar\omega}(t)$, видим, что $\varphi^\ast\bot1$ и $\varphi^\ast\in H_{\omega_{L_1}}$. Для   функции $f_\ast\in L^\psi_\beta H_{\omega_{L_1}} , $ связанной равенством (\ref{1})
с функцией $\varphi^\ast_{\bar\omega}(t),$ имеем
\begin{equation}\label{71}
 E_n(L^\psi_\beta H_{\omega_{L_1}})_{L_1}\ge E_n(f_\ast)_{L_1}
\ge
\frac{1}{\pi} \psi(n)e_n(\omega)   +O(1)\gamma_n(\psi)\omega\left(\frac 1n\right).
\end{equation}

Объединяя соотношения (\ref{37}), (\ref{39}), (\ref{70}) и (\ref{71}), получаем равенство
 (\ref{9}) для произвольных модулей непрерывности $\omega(t)$. Теорема 1 доказана.


\textbf{Contact information:}

Department of the Theory of Functions, Institute of Mathematics
of The National Academy of Sciences of Ukraine, 3, Tereshenkivska st., 01601, Kyiv, Ukraine.

E-mail: sokol@imath.kiev.ua, serdyuk@imath.kiev.ua


\begin{thebibliography}{99}

\bibitem{Stepanets2002_1}
А.И. Степанец, {\it Методы теории приближений: В
2 ч,}  Институт математики НАН Украины, Киев,  2002, Ч.\,1.


\bibitem{Stepanets1984}
  А.И. Степанец, ''Уклонения сумм Фурье на классах бесконечно дифференцируемых функций'', {\it Укр. мат. журн.,}
   \textbf{36}:6 (1984), 750--758.


\bibitem{Stepanets2001}
А.И. Степанец, ''Решение задачи
Колмогорова--Никольского для интегралов Пуассона непрерывных
функций'', {\it Матем. сб.,} {\bf192}:1 (2001), 113-138.

\bibitem{Tel1989}
С.А. Теляковский ''О приближении суммами Фурье функций высокой гладкости'',{\it Укр. мат. журн.,}  {\bf 41}:4 (1989), 510--518.

\bibitem{Falaleev2000}
 Л.П. Фалалеев ''Приближение сопряженных функций обобщенными операторами Абеля–Пуассона'', {\it Мат. заметки,}  {\bf 67}:4 (2000), 595--602.



\bibitem{Stepanets_Serdyuk_Sh_2008}
  А.И. Степанец, А.С. Сердюк, А.Л. Шидлич,   ''Классификация бесконечно дифференцируемых функций'', {\it Укр. мат. журн.,}   {\bf 60}:12 (2008), 1686--1708.

\bibitem{Korneychuk1970}
  Н.П. Корнейчук, ''Верхние грани наилучших приближений на классах дифференцируемых функций в метриках $C$ и $L$'', {\it Докл. АН СССР,} {\bf 190} (1970), 269--271.


\bibitem{Shevaldin1994}
  В.Т. Шевалдин, ''Оценки снизу поперечников классов функций, определяемых модулем непрерывности'', {\it Изв. РАН. Сер. мат.,}  {\bf58}:5 (1994),  172--188.


\bibitem{Serdyuk_Sokolenko2010_2}
  А.С. Сердюк, І.В. Соколенко, ''Лінійні методи наближення та найкращі наближення інтегралів Пуассона  функцій з класів  $H_{\omega_p}$  в метриках просторів  $L_p$'', {\it Укр. мат. журн.,}  \textbf{62}:7 (2010), 979--996.

\bibitem{Serdyuk_Sokolenko2011_1}
  A.S.~Serdyuk, I.V.~Sokolenko, ''Asymptotic behavior of best approximations of classes of Poisson integrals of functions from   $H_\omega$'', {\it Journal of Approximation Theory,} \textbf{163} (2011), 1692--1706.

\bibitem{Serdyuk_Sokolenko2011_2}
  A.S.~Serdyuk, I.V.~Sokolenko, ''Asymptotic behavior of best approximations of classes of periodic analytic functions defined by  moduli  of continuity'', {\it Proceedings of Bulgarian-Turkish-Ukrainian Scientific Conference ''Mathematical Analysis, Differential Equations and
their Applications'', Sunny Beach, Bulgaria, 15-20 September, 2010,} Sofia, Academic Publishing House ''Prof. Marin Drinov, 2011, 173--182.



\bibitem{Stepanets2002_2}
  А.И. Степанец, {\it Методы теории приближений: В
2 ч,}  Институт математики НАН Украины, Киев, 2002, Ч.\,2.



\bibitem{Serdyuk2004}
  А.С. Сердюк, ''Про один лінійний метод
наближення періодичних функцій'', {\it Проблеми теорії наближення функцій та суміжні питання: Збірник праць  Інституту математики НАН України,} Київ:   Інститут математики НАН України,  {\bf1}:1 (2004), 296--338.

\bibitem{Kornejchuk1987}
  Н.П. Корнейчук, {\it  Точные константы в теории приближения,}
 М.:  Наука, 1987.



\bibitem{Telyakovskii1964}
  С.А. Теляковский, ''Некоторые оценки для тригонометрических рядов
с квазивыпуклыми коэффициентами'', {\it Матем. сб.,} 63(105):3 (1964), 426--444.

\bibitem{Serdyuk2004_2}
  А.С. Сердюк, ''Наближення нескінченно диференційовних
пе\-рі\-о\-дич\-них функцій інтерполяційними тригонометричними
поліномами'',
{\it Укр. мат. журн.,}  \textbf{56}:4 (2004), 495--505.


\bibitem{Stechkin1980}
  С.Б. Стечкин,    ''Оценка остатка ряда Фурье для дифференцируемых функций'', {\it Тр. МИАН СССР,} 145 (1980), 126--151.

\bibitem{Telyakovskii1971}
  С.А. Теляковский, ''Оценка нормы функции через ее коэффициенты Фурье, удобная в задачах теории аппроксимации'', {\it Труды МИАН СССР,} 109 (1971), 65--97.

\end{thebibliography}
\end{document}